# An Accelerated-Decomposition Approach for Security-Constrained Unit Commitment with Corrective Network Reconfiguration— Part II: Results and Discussion

Arun Venkatesh Ramesh, *Student Member, IEEE*, Xingpeng Li, *Member, IEEE* and Kory W. Hedman, *Senior Member, IEEE*

*Abstract*— This paper presents a novel approach to handle the computational complexity in security-constrained unit commitment (SCUC) with corrective network reconfiguration (CNR) to harness the flexibility in transmission networks. This is achieved with consideration of scalability through decomposing the SCUC/SCUC-CNR formulation and then fast screening non-critical sub-problems. This is compared against the extensive formulations of SCUC and SCUC-CNR to show the advantages of the proposed typical-decomposition and accelerated-decomposition approaches to SCUC and SCUC-CNR respectively. Simulation results on the IEEE 24-bus system show that the proposed methods are substantially faster without the loss in solution quality. The proposed accelerated-decomposition approaches can be implemented for large power systems as they have great performance in the scalability tests on the IEEE 73-bus system and the Polish system when compared against the respective extensive formulations and typical-decomposition approaches. Overall, a dynamic post-contingency network can substantially alleviate network congestion and lead to a lower optimal cost.

*Index Terms*— Accelerated-decomposition approach, Benders decomposition, Corrective transmission switching, Flexible transmission, Mixed-integer linear programming, Network reconfiguration, Post-contingency congestion relief, Security-constrained unit commitment, Topology control.

## I. INTRODUCTION

Traditionally, most, if not all, industries utilize a static network ignoring network flexibility in day-ahead operations. Prior research efforts have demonstrated the benefits of harnessing the flexibility in the transmission network. Moreover, system operators practice the use of corrective network reconfiguration (CNR) to handle network congestion [1]-[2], over-voltage [3] and reliability enhancement [4]-[5] with the help of only operators' experience [6]. Prior research also pointed to economic benefits of CNR such as generation cost saving and congestion cost decrease, as well as physical benefits such as reduction of line overloads and reduced generator-startups [7] when transmission assets are treated as flexible assets and co-optimized in the day-ahead operations. Not only that, CNR can also benefit in reducing carbon emissions by alleviating congestion-induced renewable energy curtailments [8].

Network violations can be effectively handled by CNR during emergency scenarios such as a line outage. Operators can leverage an additional option where CNR removes a transmission line out of the network optimally modifying the system topology to re-route network flows to relieve post-contingency network congestion [7]. However, CNR has not been considered in day-ahead operations through security-constrained unit commitment (SCUC) owing to additional computational complexity involved with mixed-integer linear-programming constraints [9]-[15].

In Part I of this two-part paper, typical-decomposition and accelerated-decomposition approaches are proposed to address the computational complexity of both SCUC and SCUC-CNR by decomposing their extensive formulations respectively. Especially, the proposed accelerated-decomposition approach by utilizing critical sub-problem screener (CSPS) substantially reduces the solve time and outperforms the typical-decomposition approach and extensive formulation. It can effectively be scalable to handle large power systems. The proposed methods can be easily integrated into the existing practices of SCUC without the loss of solution quality in both competitive and vertical business environments. The following were the list of methods described in Part I of the paper:

- Method I: SCUC extensive formulation,
- Method II: typical-decomposition approach to SCUC,
- Method III: accelerated-decomposition approach to SCUC,
- Method IV: SCUC-CNR extensive formulation,
- Method V: typical-decomposition approach to SCUC-CNR,
- Method VI: accelerated-decomposition approach to SCUC-CNR.

Method I and Method IV serve as benchmarks for typical-decomposition and accelerated-decomposition approaches to SCUC and SCUC-CNR respectively. Method II and Method III are the proposed methods for SCUC whereas Method V and Method VI are the proposed methods for SCUC-CNR. The rest of this paper is organized as follows. Section II discusses about current industry practices in day-ahead operations in competitive markets. Section III describes the test systems considered to validate the proposed models. Following this, section IV presents the results, analysis and further discussions. Finally, section V concludes the paper.

Arun Venkatesh Ramesh and Xingpeng Li are with the Department of Electrical and Computer Engineering, University of Houston, Houston, TX, 77204, USA. Kory W. Hedman is with the School of Electrical, Computer and Energy Engineering, Arizona State University, Tempe, AZ, 85287, USA (e-mail: aramesh4@uh.edu; xingpeng.li@asu.edu; kwh@myuw.net).

## II. Industrial Day-Ahead Practices

In the United States, the wholesale energy market is a look-ahead market and consists of day-ahead and real-time markets. To ensure reliability, the North American Electric Reliability Corporation (NERC) sets several standards for the independent system operators (ISOs) to comply. Among them, the day-ahead solution must be *N*-1 compliant [16], which implies that the system solution should be capable to handle a disturbance such as a line or generator outage contingency. This is handled by committing extra generators to support the system in the case of contingencies and maintaining reserve adequacy to handle an emergency.

SCUC is performed in the day-ahead market (DAM) to obtain least-cost solution of hourly generator commitment and dispatch for the bid-cleared demand. Typical time frame of DAM covers a period of 24 hours from 00:00am to 11:59pm The SCUC requires several inputs such as load bids, generator offers, virtual bids, bilateral schedule and self-scheduling as shown in Fig. 1. In addition, the network topology and parameters are required to optimize the system for a least-cost reliable commitment and dispatch solution based on a common-pricing model. The SCUC clears almost 93%-97% of the demand, [17]-[20], following which a reliability unit commitment (RUC) is performed for meeting the forecasted loads.

Contingency analysis (CA), which is a sequence of power flow runs under different element outages, is performed by eliminating one element from the system at a time to identify any system violations for the day-ahead solution. If the CA fails, then an out of market correction is performed by committing additional generators and re-dispatching generators. The out of market correction is performed until all known violations are eliminated.

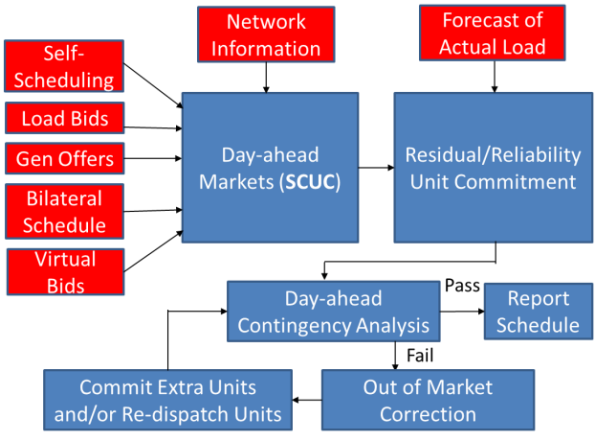

Fig. 1. ISO's typical day-ahead process.

California ISO's DAM that collects bids for energy, ancillary services, reliability unit availability, self-scheduling and virtual energy bids are open seven days prior to the operating day and closes for bids by 10:00 hours on the day prior to the operating day. Once the bids are obtained, the DAM begins with market power mitigation to identify non-competitive constraints for energy bids. Following this, the integrated forward market will clear the bids using SCUC and the RUC is used to procure additional capacity for reliability. The results for the next operational day are posted by 13:00 hours [21].

New York ISO closes the DAM for bids by 05:00 hours the day prior to the operational day. The load forecast is posted by 08:00 hours and the generator schedules are determined by clearing the energy bids and posted by 11:00 hours [22].

Midwest ISO's (MISO) DAM implements a co-optimized SCUC for energy offers and regulating reserves between 10:00-13:30 hours. MISO's DAM determines the commitment for about 1,500 resources totaling ~177,760 MW capacity, and the peak load is of ~127,125 MW [23]. After the SCUC, a rebidding is performed at 14:00 hours to run a simultaneously co-optimized security-constrained economic dispatch (SCED) for ancillary services and clearing energy prices [24].

ISO New England (ISO-NE) collects market inputs by 10:00 hours and the results are posted by 13:30 hours which publishes the generator schedules, locational marginal prices (LMP) and binding constraints. ISO-NE's network consists of 1,000+ price nodes where LMPs are calculated. The reliability of the commitment schedule is verified using a contingency analysis embedded simultaneous feasibility test to identify out-of-merit dispatches [25].

The Energy Reliability Council of Texas (ERCOT) begins DAM at 06:00 hours and ends by 18:00 hours. The information related to DAM is obtained by 06:00 hours. Then, ERCOT performs pre-market activities. The DAM clears the SCUC between 10:00-13:30 hours. Once the results for the DAM are obtained for the next operation day, the RUC begins at 14:30 hours to commit additional units by considering more accurate weather and load forecasts and updated network model. Finally, market adjustment is performed between 18:00-0:00 hours [26].

PJM's DAM collects market participant offers such as energy and regulation bids between 08:00-11:00 hours. The day-ahead results are posted by 13:30 hours after processing all the market requests from bids. After the results are available, the re-bids are processed until 14:15 hours. These re-bids and updated forecasts are used in the reliability analysis for out-of-market corrections, which goes from 14:15 hours until midnight [27].

Southwest Power Pool's (SPP's) DAM posts the available generating reserves by 06:00 hours following which SPP closes the generation offers and load bids by 09:30 hours. Between 09:30-13:00 hours the commitment and dispatch schedules are optimized using SCUC. RUC process begins at 13:45 hours after collecting re-bids. Finally, the results from RUC are posted at 16:15 hours [28].

The focus of this research is to emphasize the use of CNR to reap benefits in day-ahead operations while reducing the complexity around discrete-switching actions in the model. Though prior studies [4]-[5] point to such benefits, network reconfigurations are not extensively used in the industry owing to large system disturbances and computational complexity. However, CNR to practically handle network congestion, over-voltage and system reliability [2] have been used in real-time process. Most CNR actions are implemented in real-time with only ad-hoc procedures and operator experience as detailed by PJM operational procedures in [29]. Also protocols to enhance system reliability through eliminating internal transmission lines are presented by ISO New England [30]. Such control procedures were used during disasters like Superstorm Sandy, [31].

## III. TEST CASE DESCRIPTION

The results from the proposed Methods II-III and Methods V-VI were validated against the extensive formulation detailed in Methods I and IV, respectively, on the IEEE 24-bus system with 33 generators and 38 branches [32]. The network includes a total generation capacity of 3,393 MW and the system peak load is 2,265 MW. Furthermore, the IEEE 73-bus system and the Polish system were utilized to show the effectiveness and scalability of Methods V-VI. Table I summarizes the test systems.

The IEEE 73-bus system consists of 99 generators and 117 branches [32]. The total generation capacity is 10,215 MW and the system peak load is 8,550 MW. The Polish system is used for demonstrating the scalability of the algorithm while meeting expected performance of industry standards. It is the largest system used for this work and it consists of 2,383 buses, 327 generators and 2,895 branches [33]. The total generation capacity is 30,053 MW serving a system peak load of 21,538 MW. Two cases of the Polish system, covering a single-hour period and a 24-hour period respectively, are considered. The single-hour period case is effective to compare performance against smaller systems whereas the scalability is shown through the 24-hour period case. For the purpose of demonstrating CNR, only non-radial transmission line contingencies are considered in the *N*-1 SCUC formulation since contingency of radial lines will lead to islanding and system separation; this is consistent with industrial practice. Similarly, CNR actions, at most one action per contingency, considers only non-radial lines as possible reconfiguration actions for the same reason.

TABLE I. TEST SYSTEM SUMMARY

| System | Pgen (GW) | Pload (GW) | # bus | #gen | # branch | # radial branch |
|---|---|---|---|---|---|---|
| IEEE 24 | ~3.4 | ~2.1 | 24 | 33 | 38 | 1 |
| IEEE 73 | ~10.2 | ~8.6 | 73 | 99 | 117 | 2 |
| Polish | ~30.1 | ~21.5 | 2,383 | 327 | 2,895 | 644 |

## IV. RESULTS AND ANALYSIS

The mathematical model is implemented using AMPL and solved using Gurobi [34]-[35]. The models were run on a computer with Intel® Xeon(R) W-2195 CPU @ 2.30GHz; the CPU contains 24.75 MB of cache and 128 GB of RAM. The proposed methods were initially validated, following which sensitivity analysis, scalability and market impact are discussed.

### A. Validation for Proposed Methodologies

Since the proposed methodologies are all iterative in nature, an accuracy validation was performed to test the robustness against non-iterative extensive formulations. A MIPGAP of 0.00 was utilized on the congested network of IEEE 24-bus system for 24-hour period and the SCUC results are tabulated in Table II and SCUC-CNR results are tabulated in Table III. It was observed from Table II that the results for Method I, Method II and Method III, where CNR is not implemented, are the same. Similarly, the solutions obtained from Method IV, Method V and Method VI are the same where CNR is implemented.

The results presented in Table II and Table III prove that the proposed typical-decomposition and accelerated-decomposition methods for SCUC and SCUC-CNR are significantly faster than their respective extensive formulations for the same solution. It is intuitive that incorporating CNR into SCUC (extensive formulation without problem decomposition) will lead to additional computational complexity, which is demonstrated by the observation that the computing time of Method IV is longer than Method I. However, it is the other way for the proposed typical-decomposition and accelerated-decomposition approaches: the computational time for solving SCUC-CNR is much less than that for SCUC. The reason is that the addition of network-reconfigured post-contingency feasibility check (NR-PCFC) sub-problem in addition to the traditional post-contingency feasibility check (PCFC) in SCUC-CNR leads to increased feasibility region of the sub-problems and reduced number of cuts and iterations. When PCFC fails, the feasibility of post-contingency constraints are further verified with network reconfiguration and as a result more sub-problems are feasible when compared against Methods II-III that implement SCUC without CNR and with only PCFC. This directly translates to fewer cuts being added to the MUC after each iteration in the case of Methods V-VI, which further reduces the number of iterations before the algorithm converges. As demonstrated and discussed in detail in subsection IV.D, both Method V and Method VI of SCUC-CNR require fewer cuts and fewer iterations to converge.

TABLE II. SCUC ACCURACY ON THE IEEE 24-BUS SYSTEM

| MIPGAP=0.00 | Method I | Method II | Method III |
|---|---|---|---|
| Total cost ($) | 963,893 | 963,893 | 963,893 |
| Solve time (s) | 6,013 | 2,440 | 1351 |

TABLE III. SCUC-CNR ACCURACY ON THE IEEE 24-BUS SYSTEM

| MIPGAP=0.00 | Method IV | Method V | Method VI |
|---|---|---|---|
| Total cost ($) | 928,794 | 928,794 | 928,794 |
| Solve time (s) | 9,625 | 47 | 9 |

### B. MIPGAP Sensitivity Analysis

The MIPGAP, $\mu$, which is utilized for both the MUC and sub-problems affects the performance of all the methods. Typically, increasing the $\mu$ may lead to a less good solution, which increases the total cost. The impact of different selected $\mu$ values on total cost can be measured by the change in total cost in percentage ($\Delta Cost_\mu$), which is defined in (1).

$$\Delta Cost_\mu = \left(\frac{Cost_\mu - Cost_{\mu=0}}{Cost_{\mu=0}}\right) * 100\% \quad (1)$$

The MIPGAP sensitivity analysis is conducted on the IEEE 24-bus system and the results are shown in Fig. 2 and Fig. 3. Fig. 2 illustrates how $\Delta Cost_\mu$ varies with $\mu$ for different SCUC methods and SCUC-CNR methods. Fig. 3 shows that overall, the solve-time decreases significantly as $\mu$ increases. In general, extensive methods are more computationally intensive than the respective decomposition and accelerated-decomposition approaches except at very high relative MIPGAP such as $\mu$=1.0; higher optimality gap may result in a feasible solution for MUC faster but that solution can result in more violations in post-contingency sub-problem check which may result in additional cuts and iterations to converge. For the extensive formulations, Method IV which implements SCUC-CNR requires more computational time to execute than Method I which performs SCUC. However, the opposite is true for the typical-decomposition and the accelerated-decomposition approach; in other words, Methods V-VI that implement CNR converge faster than Methods II-III that do not implement CNR. This is because Methods V-VI, (i) requires fewer iterations and (ii) the MUC size is marginally smaller compared to Methods II-III, as explained in sub-section IV.E. Another observation is that the solve time for Method V and Method VI only reduces marginally as $\mu$ increases. One major reason is that sub-problems just verify post-contingency constraints for MUC commitment and dispatch schedule, and the solution always pointed to sub-problems being solved to $\mu=0$. Therefore, sub-problems are not affected by $\mu$. In addition, the MUC may only change marginally for the same test case.

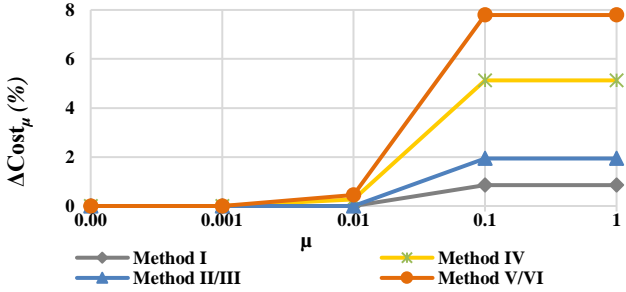

Fig. 2. $\Delta Cost_\mu$ versus relative MIPGAP $\mu$ on the IEEE 24-bus system.

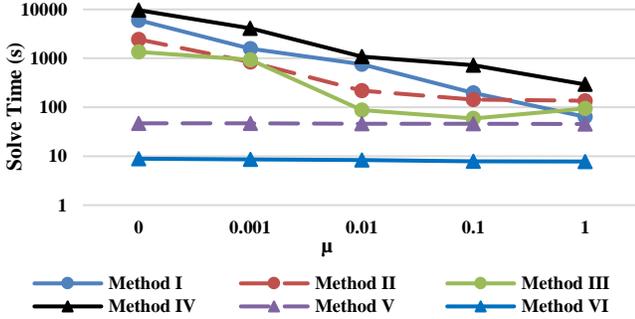

Fig. 3. Solve time versus relative MIPGAP $\mu$ on the IEEE 24-bus system.

Based on the above sensitivity analysis, $\mu=0.01$ (1%) provides a reasonable maximum cost gap of about 0.4% in a short time. For the rest of the paper, $\mu=0.01$ is used. However, it is worth noting that the performance of the proposed decomposition methodologies implementing CNR fares well under tighter tolerances if higher accuracy is required.

### C. Load Sensitivity Analysis

Four scenarios were considered: two low-load/uncongested scenarios (80%, 90%), a base-load scenario (100%) and a high-load scenario (110%). The load profile was varied using a percentage multiplied to the nodal load. Table IV shows the total cost for various methods under different load profiles.

TABLE IV. LOAD SENSITIVITY ANALYSIS ON IEEE 24-BUS SYSTEM

| Load Profile (%) | Total operational cost ($) | | | |
|---|---|---|---|---|
| | Method I | Method II/III | Method IV | Method V/VI |
| 80 | 467,883 | 467,883 | 467,883 | 467,883 |
| 90 | 624,398 | 624,398 | 623,459 | 623,459 |
| 100 | 963,893 | 963,893 | 931,919 | 932,919 |
| 110 | Infeasible | Infeasible | 1,424,140 | 1,424,140 |

In the low-load scenarios (80%, 90%), it is evident that CNR is never implemented as base-case network loading level is low and post-contingency networks are not congested. This implies all Methods I-VI obtain the same total cost.

CNR actions are observed in base-load and high-load scenarios (100%, 110%) where the network reconfiguration is utilized to relieve system congestion. This allows cheaper generators to produce more power, resulting in a reduced total operational cost. Interestingly, without CNR, the demand cannot be met due to network congestion for the high-load scenario.

### D. Intuitive Example of CNR on the IEEE 24-bus System

To explain the benefits of CNR, the following example is provided. In the IEEE 24-bus system, the line flows after the outage of line 25 are compared for SCUC and SCUC-CNR and represented in Fig. 4. It was noted that line 11 was congested in SCUC. However, in SCUC-CNR, the CNR action of disconnecting line 23 resulted in rerouting of line flows in the network that eliminated the congestion on line 11, which relieved this transmission bottleneck so that cheap power can be delivered to the demand area.

### E. Scalability Studies

One of the key research gaps is the lack of an effective algorithm for solving SCUC-CNR that is scalable for large-scale power systems and solvable in realistic time. Table V and Table VI tabulate the performance of SCUC and SCUC-CNR on IEEE 73-bus system respectively. Table V points that Method I, the extensive formulation of SCUC, requires a good starting point to solve in 7,743 seconds. One approach to have good starting solution is to utilize the commitment and dispatch results obtained from the relaxed MUC problem. However, without a starting solution, even Method I for SCUC that does not implement CNR proves to be infeasible in 100,000 seconds. A default starting solution can also be utilized where all generators are committed, which results in feasibility within 1% optimality gap in about 30,000 seconds that is still impractical. In the execution of the proposed typical-decomposition and accelerated-decomposition methods, a starting point solution is not considered. Based on Table V, the starting point has a significant influence on a large optimization problem and it can be considered for Methods II-III and Method V-VI. Since typical-decomposition approach and accelerated-decomposition approach are iterative in nature, the best starting point can be obtained from the commitment and dispatch solution in the previous iteration except for the first iteration. This may lead to further reduction in computational time.

TABLE V. SCALABILITY OF SCUC TO IEEE 73-BUS SYSTEM

| MIPGAP=0.01 | Method I | Method II | Method III |
|---|---|---|---|
| Total cost ($) | 3,224,980 | 3,223,760 | 3,223,760 |
| Solve time (s) | 7,743 | 1,273 | 367 |
| Feasibility | Feasible | Feasible | Feasible |
| Starting point | Yes | No | No |

TABLE VI. SCALABILITY SCUC-CNR TO IEEE 73-BUS SYSTEM

| MIPGAP=0.01 | Method IV | Method V | Method VI |
|---|---|---|---|
| Total cost ($) | NA | 3,218,980 | 3,218,980 |
| Solve time (s) | 100,000 | 392 | 168 |
| Feasibility | Infeasible | Feasible | Feasible |
| Starting point | Yes | No | No |

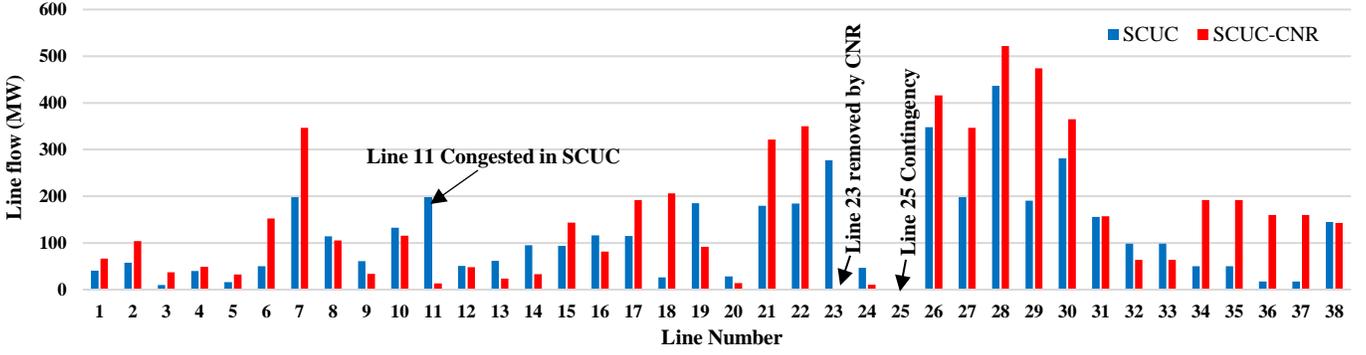

Fig. 4. Line flows in the IEEE 24-bus system under the contingency of line 25.

Table VI shows that Method IV, the extensive formulation for SCUC with CNR, lacks scalability: Method IV fails to provide a feasible solution for the IEEE 73-bus system when solved for 100,000 seconds with a good starting solution. However, this was bettered by Methods V-VI.

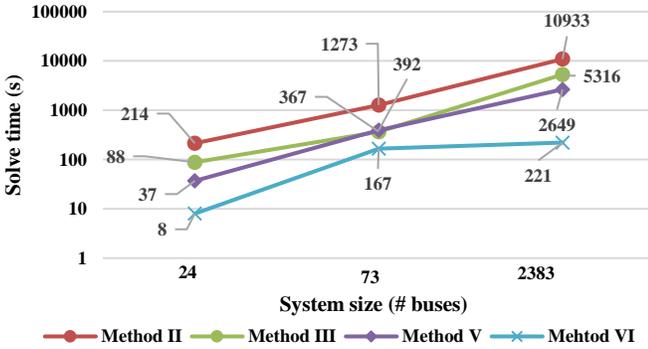

Fig. 5. Solving time versus system size.

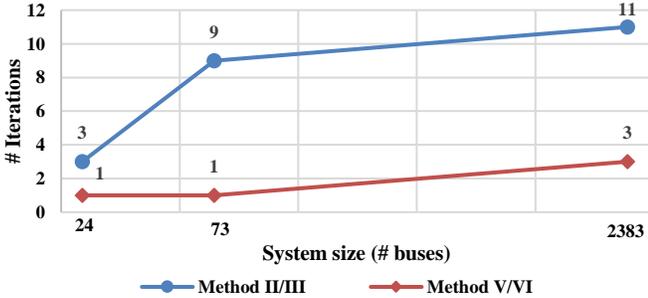

Fig. 6. Number of iterations versus size of the network.

Methods II-III and Methods V-VI are scalable to large networks such as the Polish system. Fig. 5 plots the solve time with respect to the size of the network. Methods II-III and Methods V-VI are iterative in nature and Fig. 6 plots the number of iterations to solve the problem with respect to the size of the network. Due to the size of the Polish system, the 1-hour Polish case is utilized in Fig. 5 and Fig. 6 rather than the 24-hour Polish case to compare the performance with smaller systems. Here, it is noted that Methods II-III that do not perform CNR require more iterations to converge. The transmission flexibility obtained through implementing CNR is evident from fewer iterations required to converge to a feasible solution with desired accuracy. This also means that the MUC problem that is more computationally intensive compared to sub-problems is solved fewer times, which saves a substantial amount of computational time. In addition, the number of cuts generated from infeasible post-contingency sub-problems for Methods V-VI are also less than Methods II-III. In other words, the number of constraints added to the MUC problem for each iteration for Methods V-VI is less than Methods II-III, which may lead to a less complex MUC problem and require less time to solve the MUC for each iteration for Methods V-VI. The total number of cuts added to MUC for those decomposition methods is presented in Table VII.

TABLE VII. SUB-PROBLEM AND CUT DETAILS

|        | IEEE 24-Bus | | IEEE 73-Bus | | Polish 1-Hr | |
|--------|-------------|--|-------------|--|-------------|--|
|        | Method II/III | Method IV/VI | Method II/III | Method IV/VI | Method II/III | Method IV/VI |
| # cuts | 198 | 42 | 65 | 17 | 76 | 14 |
| α      | NA  | 16 | NA | 20 | NA | 57 |

α in this table denotes number of sub-problems that were infeasible for post-contingency constraints without CNR but were feasible with CNR.

Table VIII details the results of the Polish system when it is scaled to solve for 24-hour period. Method VI, accelerated-decomposed SCUC-CNR utilizes accelerators such as the CSPS and closest branches to contingency element (CBCE), a list of 20 closest lines to the contingent line. The inclusion of accelerators in Method VI decreases the solve time by 90% as compared to Method V, decomposed SCUC-CNR, while the solution quality is retained. It is also evident that due to fewer iterations, Method VI is over 40% faster than Method III, accelerated-decomposed SCUC.

TABLE VIII. SCALABILITY TO POLISH SYSTEM FOR 24-HOUR PERIOD

| Parameters | Method III | Method V | Method VI |
|------------|-----------|----------|-----------|
| Total Cost ($) | 5,350,220 | 5,335,330 | 5,335,330 |
| € ($) | NA | 14,890 (0.28%) | 14,890 (0.28%) |
| Time (s) | 15,133.9 | 59,473.1 | 6,257.3 |
| δ | 0.04% | 0.12% | 0.12% |
| Iterations | 14 | 2 | 2 |
| # CNR | NA | 637 | 637 |
| # Cuts | 1,499 | 192 | 192 |

€ denotes the cost saving for Methods V-VI as compared to Method III. δ denotes the MIPGAP of the reported solution of MUC in the last iteration.

The comparison between Method III and Methods V-VI shows that there are 1,499 sub-problems resulting in cuts being added as constraints to the MUC problem for Method III that does not implement CNR, as opposed to 192 cuts required for Methods V-VI that implement CNR. Therefore, the MUC problem in Method III is more constrained and takes longer to solve when compared to the MUC problem in Method V-VI. Not only that, the flexibility offered by CNR is evident by the following fact: out of 829 sub-problems that failed PCFC, 637 sub-problems are feasible with CNR through NR-PCFC, which implies about 77% of contingencies that failed feasibility in post-contingency check becomes feasible when CNR actions were implemented. Moreover, Method V-VI converge faster and require only 2 iterations against Method III that requires 14 iterations, which implies the complex MUC problem is solved fewer times with Method V-VI leading to significant reduction in computational time.

The consideration of network reconfiguration for post-contingencies to alleviate network congestion in the large-scale Polish system for 24-hour period leads to a cost saving of $14,890. The discussion regarding the impact of CNR on congestion cost and markets is presented in detail in sub-section IV.F. It is to be noted that the results present an exhaustive monitoring of all non-radial transmission elements: 2,250 non-radial lines for the Polish system. This leads to 54,000 sub-problems for a 24-hour period per iteration whereas only 1,761 sub-problems were swiftly deemed as critical by CSPS. Subsequently, PCFC checked those 1,761 sub-problems and identified 829 sub-problems that failed feasibility check. The NR-PCFC that verifies feasibility of these contingencies with network reconfiguration further reduced the number of cuts required to be added to 192. Therefore, 637 sub-problems satisfied feasibility of post-contingency constraints by modifying the network topology. Though those 637 sub-problems that implemented CNR actions amounts to only 1.18% of all the sub-problems considered in the first iteration, considerable economic benefits are achieved with Methods V-VI over Method III. The solve time can be further significantly reduced if only a watch-list of key contingent lines are monitored as this will reduce the number of sub-problems drastically.

*F. Congestion Cost and Market Analysis*

The contingency-induced congestion cost, $CICC$, is calculated as the difference in total operation cost when emergency post-contingency line limits are imposed ($TC$) and not imposed ($TC_{NoEL}$) as represented in (2). The scenario when post-contingency emergency limits are not imposed is used as a benchmark since it is equivalent to implying that the system is not congested in the post-contingency situations. Method III and Method VI are considered since we are interested in calculating the amount of $CICC$ reduced when CNR is implemented.

$$CICC = TC - TC_{NoEL} \qquad (2)$$

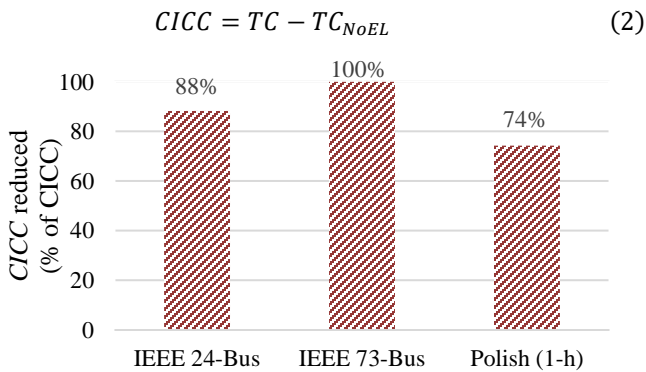

Fig. 7. *CICC* reduction with CNR for the IEEE 24-bus, IEEE 73-Bus and Polish systems.

As observed from Fig. 7, the IEEE 24-bus system was the most congested system with a contingency-induced congestion cost of $35,099 due to the considered load profile along with lower transmission capability. This was followed by the 73-bus system and 1-hour polish system with $4,550 and $ 4,150 respectively. The *CICC* is considerably reduced with CNR in all the cases by 88%, 100% and 74% respectively. This is significant in heavily congested systems as seen in the case of IEEE 24-bus system where $30,974 is saved.

TABLE IX. AVERAGE NODAL LMP ($/MWH)

| Test System | Method III | | | | Method VI | | | |
|---|---|---|---|---|---|---|---|---|
| | Mean | Min | Max | StdD | Mean | Min | Max | StdD |
| IEEE 24-Bus | 23.39 | 5.46 | 150.6 | 0.86 | 23.23 | 5.46 | 150.6 | 0.84 |
| IEEE 73-Bus | 42.75 | 9.5 | 648.4 | 1.36 | 42.19 | 4.9 | 582.4 | 1.34 |
| Polish (1-hour) | 17.72 | 15.7 | 20.8 | 0.24 | 17.56 | 17.2 | 17.8 | 0.19 |

TABLE X. LOAD PAYMENT ($)

| Test System | Method III | Method VI |
|---|---|---|
| IEEE 24 Bus | 1,171,220 | 1,112,380 |
| IEEE 73 Bus | 7,840,770 | 6,263,970 |
| Polish (1-hour) | 372,740 | 368,763 |

The market implication of reduction in *CICC* can be seen through the impact of CNR on nodal locational marginal prices (LMP). Table IX shows the average nodal LMP calculated in various systems when CNR is not used (Method III) and when CNR is implemented (Method VI). Overall it is observed that with CNR, (i) the average nodal LMP is reduced and (ii) the nodal LMP curve is flattened. It can be noted that congestion relief has a direct impact on the reduction in average nodal LMP. Similarly, it is also noted that the load payment is significantly reduced with CNR. Table X shows the total load payment for each test system with and without CNR. CNR resulted in a load payment reduction of $58,840 in the IEEE 24-bus system, $1,576,880 in the IEEE 73-bus system and $3,977 in the 1-hour Polish system, which correspond to percentage reductions of around 5.0%, 20.1% and 1.1% respectively. This makes sense since compared to the IEEE test systems, (i) the production cost of generators in the Polish system is low, (ii) the variation of system-wide generation cost in the Polish system is small, and (iii) the Polish system is lightly loaded.

*G. Forbidden Zones of Generators*

Though the forbidden zones of generators are not considered in this work, they can be integrated easily in the iterative process. Some generators may consist of a few sub-regions between the minimum and maximum outputs where the generators are unstable; those unstable sub-regions are known as forbidden zones. For instance, as shown in Fig. 8, the sub-region between $P_g^{fmin}$ and $P_g^{fmax}$ is a forbidden zone. Typically, generators are allowed to operate in these zones provided the dispatch quickly moves away from this zone. If the generators cannot operate without crossing the region, then it is operated with full ramping rates inside the zone [36].

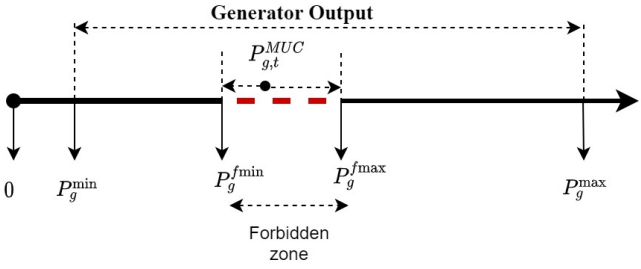

Fig. 8. Forbidden zones of operation in generators.

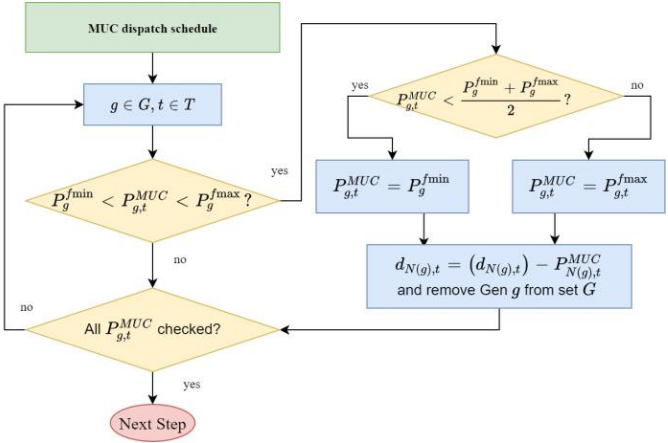

Fig. 9. Flowchart for considering forbidden zones in MUC.

Since the proposed method is an iterative process, the MUC solution after the final iteration can be verified against the forbidden zones for violations. The generator outputs for violated units can be modified as shown in Fig. 9. Once the violated generator outputs are fixed, the MUC is re-solved again for the rest of the units to satisfy the demand. Since the focus of this work is to investigate the application of CNR in SCUC, addressing generator forbidden zone is out of the scope of this paper.

### H. Renewable Energy Source Integration

The SCUC and SCUC-CNR models consider only hydro plants for this work. Other sources of variable renewable energy can be treated as traditional plants by considering a stochastic implementation which models multiple scenarios of weather (different scenarios of renewable generation profiles) as shown in [8]. This will convert deterministic SCUC and deterministic SCUC-CNR into stochastic SCUC and stochastic SCUC-CNR respectively. The commitment schedule obtained from stochastic SCUC and stochastic SCUC-CNR would be valid for all considered scenarios of weather. Since stochastic problems involve multiple scenarios, the resulting optimization problems consist of much more constraints than the associated deterministic problems. The proposed approaches to SCUC and SCUC-CNR can significantly reduce the complexity of such an extensive stochastic model by decomposing it to a MUC model that provides the base-case commitment and dispatch solution and scenario-based post-contingency constraints that can be verified as sub-problems. Though it is evident that the consideration of scenarios will lead to additional sub-problems, it was showed that the proposed typical-decomposition and accelerated-decomposition methods can handle a large number of sub-problems in a much more effective manner than directly solving a single large mixed-integer optimization problem. Note that the proposed decomposed approaches can also be used to solve stochastic SCUC and stochastic SCUC-CNR with minor modifications, which will be our future work.

### I. Other Preventive and Corrective Actions

Similar to network reconfiguration, both demand response and energy storage can increase power system flexibility and can also be used as preventive or corrective actions [37]-[40]. Though these technologies can be included explicitly by modelling constraints associated with their operational characteristics, they are not included in this paper since this paper focuses only on network reconfiguration as a corrective action. For example, the preventive action of demand response shifts demand from peak to non-peak hours in the base case, whereas as a corrective action to contingencies, non-critical demand may be curtailed for system reliability [38],[41]. Prior research, [42]-[44], states that economic benefits can be obtained as a result of using energy storage and demand response as a preventive action. However, the effect of such technologies on post-contingent scenarios are not completely studied. Especially, the effectiveness of CNR versus corrective demand response and corrective energy storage actions can be studied and their interdependence can lead to more benefits, which will be a future scope considered as an extension of this paper.

## V. CONCLUSIONS

This paper bridges the gaps by proposing an effective accelerated-decomposition approach to solving SCUC and SCUC-CNR in a competitive solve-time. An exhaustive fast screening of sub-problems was implemented, and a ranked-priority list was formed to perform CNR actions without increasing the solve-time. The proposed Method VI, accelerated-decomposed SCUC-CNR utilizing the proposed accelerators, CSPS and CBCE, can solve a large-scale power system for 24-hour period in a reasonable time. As compared to Method V that implements CNR without screening critical sub-problems, the proposed Method VI achieves a reduction of about 90% in the computational time without compromising solution accuracy.

It was noted that implementation of CNR can achieve significant cost saving and provide feasible solutions for high critical demands where there are no feasible solutions without CNR. In addition, implementation of CNR with the proposed accelerated-decomposition approach can provide quality solution much faster than the accelerated-decomposition approach without implementation of CNR since fewer iterations are required. The load payment is dramatically reduced with CNR. Load payment reduction of 1%-20% can be realized for various networks. Mainly, the advantage of the proposed accelerated-decomposition method to SCUC-CNR is that it provides quality solutions in a reasonable short time while dramatically reducing post-contingency network constraints induced congestion cost by 75%-100% in various scenarios. As a result, the total operation cost is reduced with CNR for congested networks.